# Optimal Control of Nonlinear Systems Using the Homotopy Perturbation Method

Amin Jajarmi
Department of Electrical Engineering
Ferdowsi University of Mashhad
Mashhad, Iran
jajarmi@stu-mail.um.ac.ir

Hamidreza Ramezanpour
Department of Nuclear Engineering and Physics
Amirkabir University of Technology
Tehran, Iran
h.ramezanpour@aut.ac.ir

Arman Sargolzaei, Pouyan Shafaei
Department of Electrical Engineering
Sadjad Institute of Higher Education
Mashhad, Iran
a.sargolzaei@ieee.org, pouyan.shafaei@gmail.com

*Abstract*—This paper presents a new method for solving a class of nonlinear optimal control problems with a quadratic performance index. In this method, first the original optimal control problem is transformed into a nonlinear two-point boundary value problem (TPBVP) via the Pontryagin's maximum principle. Then, using the Homotopy Perturbation Method (HPM) and introducing a convex homotopy in topologic space, the nonlinear TPBVP is transformed into a sequence of linear time-invariant TPBVP's. By solving the presented linear TPBVP sequence in a recursive manner, the optimal control law and the optimal trajectory are determined in the form of infinite series. Finally, in order to obtain an accurate enough suboptimal control law, an iterative algorithm with low computational complexity is introduced. An illustrative example demonstrates the simplicity and efficiency of proposed method.

*Index Terms*—nonlinear optimal control problem, Pontryagin's maximum principle, two-point boundary value problem, Homotopy Perturbation Method

## I. Introduction

Theory and application of optimal control has been widely used in different fields such as biomedicine [1], aircraft systems [2], robotic [3], etc. However, optimal control of nonlinear systems is a challenging task which has been studied extensively for decades.

Methods of solving nonlinear optimal control problems (OCP's) can be divided into two categories. The first category, which contains direct methods, converts the problem into a nonlinear programming by using the discretization or parameterization techniques [4-5]. The second category contains indirect methods and leads to the Hamilton-Jacobi-Bellman (HJB) equation, based on dynamic programming [6], or nonlinear two-point boundary value problem (TPBVP), based on the Pontryagin's maximum principle [7]. In general, the HJB equation is a nonlinear partial differential equation that is hard to solve in most cases. An excellent literature review on the methods for approximating the solution of HJB equation is provided in [8]. Besides, nonlinear TPBVP has no analytical solution except for a few simple cases. Thus, many researches have been devoted to find an approximate solution for the nonlinear TPBVP's. Recently, successive approximation approach (SAA) and sensitivity approach have been introduced in [9] and [10], respectively. In those, a sequence of nonhomogeneous linear time-varying TPBVP's is solved instead of directly solving the nonlinear TPBVP derived from the Pontryagin's maximum principle. However, solving time-varying equations is much more difficult than solving time-invariant ones.

The Homotopy Perturbation Method (HPM) was initially proposed by the Chinese mathematician J. H. He [11-12]. This method has been widely used to solve nonlinear problems in different fields [13-15]. In contrast to the perturbation method [16], the HPM is independent upon small/large physical parameters in system model. However, like the other traditional non-perturbation methods such as the Lyapunov's artificial small parameter method [17] and Adomian's decomposition method [18], uniformly convergence of the solution series obtained via the HPM can not be ensured.

In this paper, based on the HPM, a new method is proposed to solve a class of nonlinear OCP's. In this method, first the original nonlinear OCP is transformed into a nonlinear TPBVP by using the maximum principle. Applying the HPM transforms the nonlinear TPBVP into a sequence of linear time-invariant TPBVP's. Solving the proposed linear TPBVP sequence in a recursive manner leads to the optimal control law and the optimal trajectory in the form of infinite series. The proposed method avoids the trouble of directly solving the nonlinear TPBVP or the HJB equation. In addition, it avoids solving a sequence of linear time-varying TPBVP's. It only requires solving a sequence of linear time-invariant TPBVP's. Hence, it reduces the computational complexity, effectively. The rest of paper is organized as follows. In section 2 the statement of problem is discussed. Section 3 elaborates the proposed method. In order to obtain an accurate enough suboptimal control law, an efficient algorithm with low computational complexity is introduced in section 4. Section 5 contains a numerical example to show the effectiveness of proposed method. Finally, conclusions and future works are given in the last section.

## II. Problem Statement

Consider the following nonlinear control system:

$$\begin{cases} \dot{x}(t) = Ax(t) + Bu(t) + f(x(t)) \\ x(t_0) = x_0, x(t_f) = x_f \end{cases} \quad (1)$$

where $A$ and $B$ are real constant matrices of appropriate dimensions, $x \in R^n$ is the state vector and $u \in R^m$ is the control vector, $f$ is a nonlinear polynomial vector function where $f(0) = 0$, $x_0 \in R^n$ and $x_f \in R^n$ are the initial and final states, respectively. The objective is to find the optimal control law $u^*(t)$ that minimizes the following quadratic performance index:

$$J = \frac{1}{2}\int_{t_0}^{t_f} \left( x^T(t) Q x(t) + u^T(t) R u(t) \right) dt \quad (2)$$

subject to the system (1) where $Q \in R^{n \times n}$ and $R \in R^{m \times m}$ are positive semi-definite and positive definite matrices, respectively.

According to the Pontryagin's maximum principle, the optimality condition is obtained as the following nonlinear TPBVP:

$$\begin{cases} \dot{x}(t) = Ax(t) - BR^{-1}B^T \lambda(t) + f(x(t)) \\ \dot{\lambda}(t) = -Qx(t) - A^T \lambda(t) - f_x(x(t))\lambda(t) \\ x(t_0) = x_0, x(t_f) = x_f \end{cases} \quad (3)$$

where $f_x = \dfrac{\partial f}{\partial x}$ and $\lambda \in R^n$ is the co-state vector. Also, the optimal control law is given by:

$$u^*(t) = -R^{-1}B^T \lambda(t). \quad (4)$$

## III. Proposed Method

Unfortunately, solving the nonlinear TPBVP (3) is very difficult in most cases. In order to overcome this difficulty, in this section, we introduce a new method, based on the HPM, which transforms the nonlinear TPBVP (3) into a sequence of linear time-invariant TPBVP's.

Let define the operators $F_1(x(t), \lambda(t))$ and $F_2(x(t), \lambda(t))$ as follows:

$$F_1(x(t), \lambda(t)) \stackrel{\Delta}{=} \dot{x}(t) - Ax(t) + BR^{-1}B^T \lambda(t) - f(x(t)) \quad (5)$$

$$F_2(x(t), \lambda(t)) \stackrel{\Delta}{=} \dot{\lambda}(t) + Qx(t) + A^T \lambda(t) + f_x(x(t))\lambda(t) \quad (6)$$

From (3) it is obvious that:

$$F_i(x(t), \lambda(t)) = 0 : i = 1,2 \quad (7)$$

The operator $F_i : i = 1,2$ can generally be divided into two parts, a linear part and a nonlinear part. So, we can write:

$$F_i(x(t), \lambda(t)) = L_i(x(t), \lambda(t)) + N_i(x(t), \lambda(t)) : i = 1,2 \quad (8)$$

where $L_i$ and $N_i$ are respectively the linear and nonlinear parts of $F_i$ for $i = 1,2$. Now, we construct a homotopy for (8) as follows:

$$\begin{cases} (1-p)L_1(\tilde{x}(t,p), \tilde{\lambda}(t,p)) + pF_1(\tilde{x}(t,p), \tilde{\lambda}(t,p)) = 0 \\ (1-p)L_2(\tilde{x}(t,p), \tilde{\lambda}(t,p)) + pF_2(\tilde{x}(t,p), \tilde{\lambda}(t,p)) = 0 \end{cases} \quad (9)$$

with boundary conditions:

$$\tilde{x}(t_0, p) = x_0, \tilde{x}(t_f, p) = x_f \quad (10)$$

where $p \in [0,1]$ is an embedding parameter which is called homotopy parameter. Setting $p = 0$ and $p = 1$ in (9) yields:

$$p = 0 \Rightarrow L_i(\tilde{x}(t,0), \tilde{\lambda}(t,0)) = 0 : i = 1,2 \quad (11)$$

$$p = 1 \Rightarrow F_i(\tilde{x}(t,1), \tilde{\lambda}(t,1)) = 0 : i = 1,2 \quad (12)$$

Therefore, if the homotopy parameter $p$ changes from zero to unity, $\tilde{x}(t, p)$ and $\tilde{\lambda}(t, p)$ change from the solution of (11) to the solution of (12). In topology we call it deformation. Obviously, when $p = 1$, TPBVP (9)-(10) is equivalent to the nonlinear TPBVP (3).

**Theorem 3.1.** The solution of nonlinear TPBVP (3) can be written as $x(t) = \sum_{n=0}^{\infty} x^{(n)}(t)$ and $\lambda(t) = \sum_{n=0}^{\infty} \lambda^{(n)}(t)$ where the $n$-th order terms $x^{(n)}(t)$ and $\lambda^{(n)}(t)$ for $n \geq 0$ are achieved recursively by solving a sequence of linear time-invariant TPBVP's.

**Proof.** Assume that the embedding parameter $p$ is a small parameter and $\tilde{x}(t, p)$ and $\tilde{\lambda}(t, p)$ are infinitely differentiable with respect to $p$ around $p = 0$. Expanding $\tilde{x}(t, p)$ and $\tilde{\lambda}(t, p)$ as Maclaurin series yields:

$$\begin{cases} \tilde{x}(t, p) = x^{(0)}(t) + x^{(1)}(t)p + x^{(2)}(t)p^2 + \cdots \\ \tilde{\lambda}(t, p) = \lambda^{(0)}(t) + \lambda^{(1)}(t)p + \lambda^{(2)}(t)p^2 + \cdots \end{cases} \quad (13)$$

where $x^{(n)}(t) = \dfrac{1}{n!} \dfrac{\partial^n \tilde{x}(t,p)}{\partial p^n}\bigg|_{p=0}$ and

$\lambda^{(n)}(t) = \dfrac{1}{n!} \dfrac{\partial^n \tilde{\lambda}(t,p)}{\partial p^n}\bigg|_{p=0}$. Substituting (13) in (9) we obtain:

$$\begin{cases} L_i(x^{(0)}(t), \lambda^{(0)}(t)) \\ + p\big(L_i(x^{(1)}(t), \lambda^{(1)}(t)) + N_i(x^{(0)}(t), \lambda^{(0)}(t))\big) \\ + p^2\big(L_i(x^{(2)}(t), \lambda^{(2)}(t)) + N_i(x^{(1)}(t), \lambda^{(1)}(t))\big) \\ + \cdots = 0 \end{cases} \quad (14)$$

From (14) we can easily obtain:

$$\begin{cases} L_1(x^{(0)}(t), \lambda^{(0)}(t)) = 0 \\ L_2(x^{(0)}(t), \lambda^{(0)}(t)) = 0 \\ x^{(0)}(t_0) = x_0, x^{(0)}(t_f) = x_f \end{cases} \quad (15a)$$

$$\begin{cases} L_1(x^{(n)}(t), \lambda^{(n)}(t)) + N_1(x^{(n-1)}(t), \lambda^{(n-1)}(t)) = 0 \\ L_2(x^{(n)}(t), \lambda^{(n)}(t)) + N_2(x^{(n-1)}(t), \lambda^{(n-1)}(t)) = 0 \\ x^{(n)}(t_0) = 0, x^{(n)}(t_f) = 0 \\ n \geq 1 \end{cases} \quad (15b)$$

Therefore, at each step, a nonhomogeneous linear time-invariant TPBVP is obtained in which nonhomogeneous terms are calculated using the information obtained from previous step. Consequently, the original nonlinear TPBVP (3) has been transformed into a sequence of linear time-invariant TPBVP's which should be solved in a recursive process.

After obtaining $x^{(n)}(t)$ and $\lambda^{(n)}(t)$ for $n \geq 0$, we should set $p=1$ in (13) to obtain the exact solution of problem (3). Setting $p=1$ in (13) yields:

$$\begin{cases} x(t) = \tilde{x}(t,1) = x^{(0)}(t) + x^{(1)}(t) + x^{(2)}(t) + \cdots \\ \lambda(t) = \tilde{\lambda}(t,1) = \lambda^{(0)}(t) + \lambda^{(1)}(t) + \lambda^{(2)}(t) + \cdots \end{cases} \quad (16)$$

and the proof is complete.

**Remark 3.1.** It should be noted that series in (16) converge rapidly for most cases; however, convergence rate depends upon the nonlinear operators [11].

**Remark 3.2.** Substituting (16) in (4), the optimal control law is obtained as follows:

$$u^*(t) = -R^{-1}B^T \sum_{n=0}^{\infty} \lambda^{(n)}(t). \quad (17)$$

## IV. SUBOPTIMAL CONTROL DESIGN

In fact, obtaining the optimal control law as in (17) is almost impossible since (17) contains infinite series. Therefore, in practical applications, by replacing $\infty$ with a finite positive integer $M$ in (17), an $M$-th order suboptimal control law is obtained as follows:

$$u_M(t) = -R^{-1}B^T \sum_{n=0}^{M} \lambda^{(n)}(t) \quad (18)$$

The integer $M$ in (18) is generally determined according to a concrete control precision. For example, every time $x^{(n)}(t)$ and $\lambda^{(n)}(t)$ are obtained from the presented linear TPBVP sequence in (15a)-(15b), we let $M = n$ and calculate the $M$-th order suboptimal control law from (18). Then, the following quadratic performance index can be calculated:

$$J^{(M)} = \frac{1}{2} \int_{t_0}^{t_f} \left( x^T(t)Qx(t) + u_M^T(t)Ru_M(t) \right) dt \quad (19)$$

where $u_M(t)$ has been obtained from (18) and $x(t)$ is the corresponding state trajectory obtained from applying $u_M(t)$ to the original nonlinear system in (1) with $x(t_0) = x_0$. The $M$-th order suboptimal control law has desirable accuracy if for given positive constant $\varepsilon > 0$, the following condition holds:

$$\left| J^{(M)} - J^{(M-1)} \right| < \varepsilon \quad (20)$$

In order to obtain an accurate enough suboptimal control law, we present an iterative algorithm with low computational complexity as follows:

**Algorithm:**
**Step 1.** Construct a homotopy as (9).
**Step 2.** Let $n = 0$.
**Step 3.** Calculate the $n$-th order terms $x^{(n)}(t)$ and $\lambda^{(n)}(t)$ from the presented linear TPBVP sequence in (15a)-(15b).
**Step 4.** Let $M = n$ and obtain the $M$-th order suboptimal control law $u_M(t)$ from (18), apply it to the original nonlinear system with $x(t_0) = x_0$ to obtain the corresponding state trajectory $x(t)$, and then calculate $J^{(M)}$ according to (19).
**Step 5.** If (20) holds for the given small enough constant $\varepsilon > 0$, go to step 6; else replace $n$ by $n+1$ and go to step 3.
**Step 6.** Stop the algorithm; $u_M(t)$ is the desirable suboptimal control law.

## V. NUMERICAL EXAMPLE

In this example, the optimal maneuvers of a rigid asymmetric spacecraft are considered as the following nonlinear OCP:

Min $J = \frac{1}{2} \int_0^{100} \left( u_1^2 + u_2^2 + u_3^2 \right) dt$

s.t.

$$\begin{cases} \dot{\omega}_1(t) = -\frac{(I_3 - I_2)}{I_1} \omega_2 \omega_3 + \frac{u_1}{I_1}, \omega_1(0) = 0.01 \text{ r/s} \\ \dot{\omega}_2(t) = -\frac{(I_1 - I_3)}{I_2} \omega_1 \omega_3 + \frac{u_2}{I_2}, \omega_2(0) = 0.005 \text{ r/s} \\ \dot{\omega}_3(t) = -\frac{(I_2 - I_1)}{I_3} \omega_1 \omega_2 + \frac{u_3}{I_3}, \omega_3(0) = 0.001 \text{ r/s} \\ \omega_1(100) = \omega_2(100) = \omega_3(100) = 0 \text{ r/s} \end{cases} \quad (21)$$

where $\omega_1$, $\omega_2$, and $\omega_3$ are angular velocities of the spacecraft, $u_1$, $u_2$, and $u_3$ are control torques, $I_1 = 86.24 \text{ kg m}^2$, $I_2 = 85.07 \text{ kg m}^2$, $I_3 = 113.59 \text{ kg m}^2$ are the spacecraft principle inertia.

In order to obtain an accurate enough suboptimal control law, we applied the proposed algorithm with tolerance error bound $\varepsilon = 10^{-12}$. In this case, convergence has been achieved after 4 iterations, i.e. $\left|J^{(4)} - J^{(3)}\right| = 5.77 \times 10^{-13} < 10^{-12}$, and a minimum of $J^{(4)} = 0.004689$ has been obtained. Simulation results are shown in figures 1 and 2 as compared with those obtained by collocation method. As it is seen, results of both methods are very close to each other. This confirms that the proposed method yields excellent results. Furthermore, in contrast to the collocation method, computing procedure of our method is very straightforward that can be done by pencil-and-paper only.

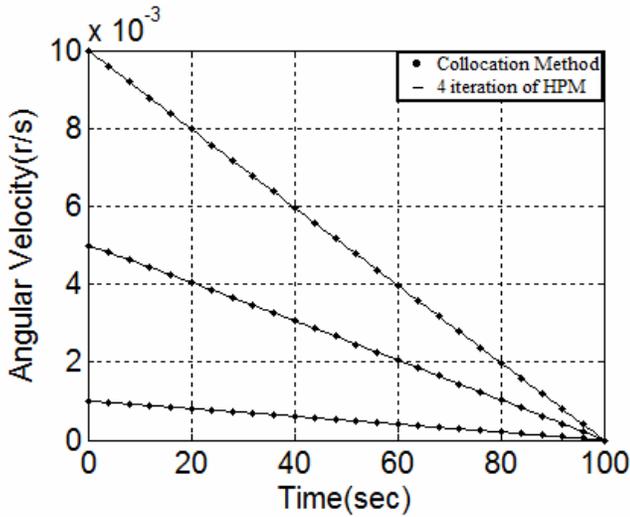

Figure 1.  Suboptimal state trajectories

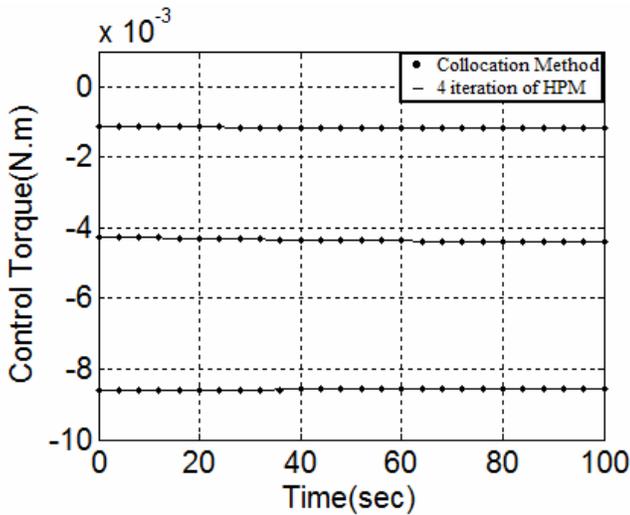

Figure 2.  Suboptimal control laws

## VI. CONCLUSION

In this paper, based on the HPM, an efficient method has been introduced to solve a class of nonlinear OCP's. In this method, by introducing a recursive process, the optimal control law is determined in the form of infinite series with easy-computable terms. The proposed method avoids directly solving the nonlinear TPBVP or the HJB equation. In addition, despite of the successive approximation approach [9] and sensitivity approach [10], it avoids solving a sequence of linear time-varying TPBVP's. It only requires solving a sequence of linear time-invariant TPBVP's. Therefore, in view of computational complexity, the proposed method is more practical than the above-mentioned approximate methods.

Future works are focused on extending the method for solving more general form of nonlinear OCP.

REFERENCES

[1] M. Itik, M. U. Salamci, and S. P. Banksa, "Optimal control of drug therapy in cancer treatment," Nonlinear Analysis, vol. 71, pp. 473-486, 2009.

[2] W. L. Garrard and J. M. Jordan, "Design of nonlinear automatic flight control systems," Automatica, vol. 13, no. 5, pp. 497-505, 1977.

[3] S. Wei, M. Zefran, and R. A. DeCarlo, "Optimal control of robotic system with logical constraints: application to UAV path planning," IEEE International Conference on Robotic and Automation, Pasadena, Ca, Usa, 2008.

[4] O. Stryk and R. Bulirsch, "Direct and indirect methods for trajectory optimization," Ann. Oper. Res., vol. 37, pp. 357-373, 1992.

[5] C. J. Goh and K. L. Teo, "Control parameterization: a unified approach to optimal control problem with general constraints," Automatica, vol. 24, pp. 3-18, 1988.

[6] R. Bellman, On the theory of dynamic programming, P. Natl. Acad. Sci. USA, vol. 38, no. 8, pp. 716-719, 1952.

[7] L. S. Pontryagin, Optimal control processes, Usp. Mat. Nauk, vol. 14, pp. 3-20, 1959.

[8] R. W. Beard, G. N. Saridis, and J. T. Wen, Galerkin approximations of the generalized Hamilton-Jacobi-Bellman equation, Automatica, vol. 33, no. 12, pp. 2159-2177, 1997.

[9] G. Y. Tang, "Suboptimal control for nonlinear systems: a successive approximation approach," Syst. Control Lett., vol. 54, no. 5, pp. 429-434, 2005.

[10] G. Y. Tang, H. P. Qu, and Y. M. Gao, "Sensitivity approach of suboptimal control for a class of nonlinear systems," J. Ocean Univ. Qingdao, vol. 32, no. 4, pp. 615-620, 2002.

[11] J. H. He, "Homotopy perturbation technique," Computational Methods in Applied Mechanics and Engineering, vol. 178, pp. 257–262, 1999.

[12] J. H. He, "A coupling method of a homotopy technique and a perturbation technique for non-linear problems," International Journal of Non-Linear Mechanics, vol. 35, pp. 37–43, 2000.

[13] J. H. He, "Homotopy perturbation method for solving boundary value problems," Phys. Lett. A, vol. 350, pp. 87–88, 2006.

[14] D. D. Ganji and A. Sadighi, "Application of He's homotopy perturbation method to nonlinear coupled systems of reaction–diffusion equations," Int. J. Nonlinear Sci. Numer. Simul., vol. 7, pp. 411–418, 2006.

[15] M. Rafei and D. D. Ganji, "Explicit solutions of Helmholtz equation and fifth-order KdV equation using homotopy perturbation method," Int. J. Nonlinear Sci. Numer. Simul., vol. 7, pp. 321–328, 2006.

[16] J. A. Murdock, Perturbations: Theory and Methods, Classics in Applied Mathematics, vol. 27, SIAM, 1999.

[17] A. M. Lyapunov, General Problem on Stability of Motion, Taylor & Francis, London, 1992 (English translation).

[18] G. Adomian, and G. E. Adomian, A global method for solution of complex systems, Math. Modelling, vol. 5, no. 4, pp. 251-263, 1984.